\documentclass[11pt,a4paper]{article}

  \usepackage[a4paper,text={450pt,650pt},centering]{geometry}

\usepackage{amsmath,amsthm} 
\usepackage{amsfonts} 
\usepackage{amssymb}
\usepackage{graphicx} 

\usepackage[utf8]{inputenc}
\usepackage[english]{babel}

\usepackage[disable,textsize=tiny,textwidth=1.7cm]{todonotes}

\allowdisplaybreaks[3]

\numberwithin{equation}{section}

\usepackage[font=small,labelfont=bf,width=0.85\textwidth]{caption}

\DeclareMathSymbol{\Gamma}{\mathalpha}{letters}{"00}
\DeclareMathSymbol{\Delta}{\mathalpha}{letters}{"01}
\DeclareMathSymbol{\Theta}{\mathalpha}{letters}{"02}
\DeclareMathSymbol{\Lambda}{\mathalpha}{letters}{"03}
\DeclareMathSymbol{\Xi}{\mathalpha}{letters}{"04}
\DeclareMathSymbol{\Pi}{\mathalpha}{letters}{"05}
\DeclareMathSymbol{\Sigma}{\mathalpha}{letters}{"06}
\DeclareMathSymbol{\Upsilon}{\mathalpha}{letters}{"07}
\DeclareMathSymbol{\Phi}{\mathalpha}{letters}{"08}
\DeclareMathSymbol{\Psi}{\mathalpha}{letters}{"09}
\DeclareMathSymbol{\Omega}{\mathalpha}{letters}{"0A}

\ifx\genfrac\sdflkaj\else\fi

\newcommand{\set}[1]{\{#1\}}

\newcommand{\hypref}[2]{\ifx\href\asklfhas #2\else\href{#1}{#2}\fi}

\newcommand{\figref}[1]{Fig.~\ref{#1}}

\newcommand{\thmref}[1]{Theorem~\ref{#1}}

\ifx\href\asklfhas\newcommand{\href}[2]{#2}\fi

\newtheorem{theorem}{Theorem}

\newtheorem{prop}{Proposition}

\newtheorem{definition}{Definition}

\newcommand{\st}{\, :\,}

\newcommand{\kre}{Kreweras}
\newcommand{\gess}{Gessel}
\newcommand{\rkre}{Reverse Kreweras}
\newcommand{\rk}{$\bar{\text{K}}$}
\newcommand{\rkm}{\bar{\text{K}}}
\newcommand{\mr}{Mishna-Rechnitzer}

\newcommand{\stepset}{\mathcal{S}}

\newcommand{\xb}{\bar{x}}
\newcommand{\yb}{\bar{y}}

\newcommand{\integer}{\mathbb{Z}}

\newcommand{\inner}[2]{\langle #1,#2\rangle}
\newcommand{\rootsys}{R}
\newcommand{\rgen}{\sigma}
\newcommand{\chamber}{\mathcal{C}}
\newcommand{\len}{\ell}
\newcommand{\cpath}{\text{Walk}}
\newcommand{\fpath}{\text{Free}}

\newcommand{\cwall}{\bar\chamber}

\newcommand{\group}{\mathcal{G}}

\newcommand{\lattice}{\mathcal{L}}
\newcommand{\qpln}{\mathcal{Q}}
\newcommand{\boundary}{\mathcal{Q}_\partial}

\newcommand{\ba}{\mathbf{a}}
\newcommand{\bb}{\mathbf{b}}

\newcommand{\gstepcount}{\mathcal{Z}^\text{P}}

\newcommand{\gfreecount}{\mathcal{F}^\text{P}}

\newcommand{\bm}{Bousquet-Mélou}

\newcommand{\act}{\boldsymbol{\cdot}}
\newcommand{\eqdef}{\mathrel{\mathop :}=}

\begin{document}

\title{Weyl chambers for short step  Quarter-plane Lattice Paths.}
\author{Richard Brak}

\maketitle

\begin{center}
    

\textit{School of Mathematics \\%
The University of Melbourne\\%
Parkville, VIC 3052,\\ Australia}\vspace{3mm}%

\verb+rb1@unimelb.edu.au+ %
 
\par\vspace{1cm}

\begin{abstract} We consider four examples of short step lattice paths confined to the quarter plane. These are the \kre,  \rkre, \gess, and \mr\ lattice paths.   
The \rkre\ are straightforward to solve and thus interesting as a contrast to the \kre\  paths and \gess\ paths   as the latter two have historically been \emph{significantly} more difficult to solve. The \mr\ paths are interesting as they are associated with and infinite order  group.   
We will give some geometrical insight into all these properties by considering the Weyl chambers associated with their step sets. 

 For \rkre\ paths the   Weyl chamber walls coincide  with the  quarter plane boundary and  hence the problem is readily solvable by Bethe  Ansatz or by using the Gessel-Zeilberger Theorem. 
 For \kre\ paths the quarter plane corresponds to the union of two adjacent Weyl Chambers and hence neither the Bethe Ansatz nor the    Gessel-Zeilberger Theorem are directly applicable making the problem considerably more difficult  to solve. Similarly, the quarter plane for   \gess\ paths is the union of three Weyl chambers.
 
 For \mr\ paths the step set has non-zero barycenter leading to an affine dihedral reflection group. The affine structure corresponds to the drift in the random walk. The quarter plane is the union of an infinite number of Weyl alcoves.
\end{abstract}
\end{center}
 \vfill

\newpage
\renewcommand{\thefootnote}{\arabic{footnote}}
\setcounter{footnote}{0}


\section{Introduction and definitions}


The  lattice $\lattice$ is  $\lattice=\integer\times\integer$ where $\integer$ is the set of integers.  
The \textbf{quarter plane}, $\qpln$, is the subset $\qpln =\set{(a,b)\in \lattice \st a>0\,\, \text{and}\,\, b>0}$. 
Define the $y$-\textbf{boundary} to be the subset $\qpln_y=\set{(0,b)\in \lattice \st  b\ge0  }$ and the $x$-\textbf{boundary}  to be  the subset $\qpln_x=\set{(a,0)\in \lattice \st  
a\ge0  }$.  
The \textbf{boundary} of the quarter plane, $\boundary$ is defined as the subset of points $\boundary=\qpln_x \cup \qpln_y $.

A \textbf{step set}, $\stepset$, is any finite subset of $\lattice$. 
A \textbf{step} is an element of $\stepset$. A \textbf{short step}  is a step $(a,b)$ for which $a^2+b^2\le2$ ie.\ a nearest-neighbour step. 
A \textbf{lattice path} or \textbf{random walk} of length $t$, or $t$-path,  with step set $\stepset$  is a finite sequence $s_1s_2s_3\dots s_t$ of steps, $s_i\in\stepset$. 
If a path starts at $\ba\in \lattice$ and has steps   $s_1s_2s_3\dots s_t$ then it is said to end at $\bb=\ba+ \sum_{i=1}^t s_i$. The \textbf{step set generating function} is the Laurent polynomial
\begin{equation}
    \lambda(\stepset)=\sum_{(a,b)\in\stepset}x^a y^b\,.
\end{equation}

We wish to compute the number of length $t$ lattice paths that start at some point $\ba=(a_1,a_2)\in \qpln$ and end at some point $\bb=(b_1,b_2)\in \qpln$ with the condition that the path is not allowed to step on the quarter plane boundary ie.\ \emph{all} of the partial sums  $\ba+ \sum_{i=1}^p s_i$, $1\le p\le t$ are points in $\qpln$. 
We will refer to such paths as \textbf{quarter plane paths} \cite{Fayolle:1999qy}.
Note, this definition differs slightly from that used in \cite{Mishna:2009kx} where the quarter plane boundary is shifted off the $x$ and $y$ axes.

Let
 \begin{equation}
     \gstepcount_t(\ba\to \bb)=
     \text{Number of  $t$-paths  in $\qpln$ from $\ba$ to $\bb$, with step set $\stepset_P$}
 \end{equation}
We will also consider paths  in $\lattice$    \emph{without} the quarter plane boundary constraint. We will refer to these a \textbf{free  lattice paths} with
 \begin{equation}
     \gfreecount_t(\ba\to \bb)=
     \text{Number of free $t$-paths in $\lattice$ from $\ba$ to $\bb$, with step set $\stepset_P$}
 \end{equation}
Note, for free paths $\ba$ and $\bb$ are no longer constrained to the quarter plane.
We define the $t$-step \textbf{free generating function}
\begin{equation} 
    G_t^\text{P}(\ba;x,y)=\sum_{\bb\in\lattice} \gfreecount_t(\ba\to \bb)x^{b_1} y^{b_2}
\end{equation}
where $\bb=(b_1,b_2)$. 


\bm\ and Misha \cite{Mishna:2009kx}   considered  56 different short step quarter plane  path problems. Integral to the many solutions obtained by \bm\ and Misha was the use of a  ``substituion'' group associated with the step set. This group is defined as follows: Let $g_1$ and $g_2$ be rational functions of $x$ and $y$.
The action of the  map \mbox{$g: (x,y)\mapsto (g_1(x,y),g_2(x,y))$},   on $\lambda(\stepset)$ is defined as $g\act \lambda(\stepset) \eqdef \sum_{(a,b)\in\stepset}g_1(x,y)^a g_2(x,y)^b$. The set of such maps that leave $\lambda$ invariant i.e.\ $ g\act  \lambda  = \lambda$, form  a group (under composition).

In this paper we will revisit four of the 56  cases. In particular, the three cases are \kre(K) \cite{Kreweras:1965aa},  \rkre(\rk), \gess(G), and \mr(MR) paths \cite{Mishna:2009kx} which have respective step sets (illustrated below):
\begin{align}
    \text{\rkre\ :} &\qquad \stepset_{\rkm}=\set{ (1,0),(0,1),(-1,-1)},\\
    \text{\kre\ :} &\qquad \stepset_\text{K}=\set{ (-1,0),(0,-1),(1,1)} \\
      \text{\gess\ :} &\qquad \stepset_\text{G}=\set{(-1,0),(1,0), (-1,-1),(1,1)} \\
    \text{\mr\ :} &\qquad \stepset_{\text{MR}}=\set{ (-1,1),(1,-1),(1,1)} 
\end{align}
\begin{center}
    \includegraphics[width=10cm]{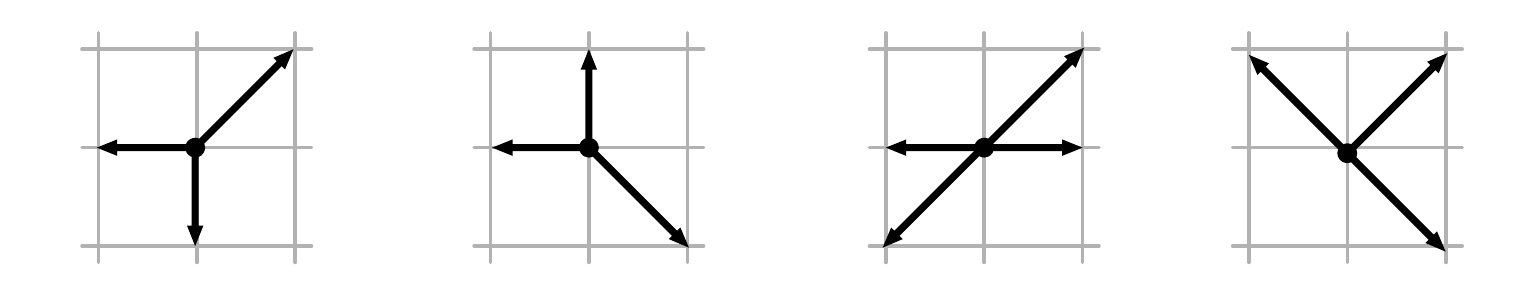}
\end{center}
and respective step set generating functions:
\begin{align}
  \text{\rkre\ :} &\qquad \lambda_{\rkm}(x,y)   = x+y+\xb\yb  \\
    \text{\kre\ :} &\qquad  \lambda_K(x,y)   = \xb+\yb+xy  \\
        \text{\gess\ :} &\qquad  \lambda_G(x,y)   =x+ \xb+\xb \yb+xy  \\
    \text{\mr\ :} &\qquad  \lambda_{MR}(x,y)   = \xb y + x \yb+xy
\end{align}
where we use the notation
\begin{equation}
    \xb \eqdef  \frac{1}{x},\qquad \yb \eqdef \frac{1}{y}\,.
\end{equation}

\kre\  paths and \gess\ paths are   interesting as they have historically been significantly more difficult to solve exactly \cite{bousquet-melou:2005uq}. 
We will give some geometrical insight into why this is the case. 

The \mr\ paths are of interest as \rm \cite{Mishna:2009kx} and \cite{bousquet-melou:2005uq} found an infinite substitution group associated with the step set but did not characterise it.


In this paper, rather than substitution groups, we consider   reflection groups (acting on $\lattice$ rather than $\lambda$) associated with the various step sets.  All the finite reflection groups turn out to be finite dihedral reflection groups (isomorphic to the the corresponding substitution group).

Of particular relevance are the Weyl chambers associated with each group.    
We shall see that the relationship between the Weyl chamber of the refection group and the quarter plane is integral to the method of solution. 
In particular if the Weyl chamber corresponds to the quarter plane then the Bethe Ansatz  method \cite{brak:1998un} and Gessel-Zeilberger theorem \cite{gessel:1992kx} are readily applied. 
This is the case for the \rkre\ paths but \emph{not} for  \kre\ or \gess\  paths. 

How the infinite reflection group   associated with the  \mr\ step set can be used to solve the quarter plane problem is not clear. 
In this context the Bethe Ansatz manifests itself as the same group orbit sum \eqref{eq_orbit} that appears in the Gessel-Zeilberger theorem \cite{gessel:1992kx}.

We briefly summarise the Gessel-Zeilbger theorem. 
Let $\rootsys$ be a finite or affine root system, let $\group$ be a Weyl group and $\Delta$ any of basis for $\rootsys$. 
The length of an  element   $g\in \group$, denoted $\len(g)$, is the least number of terms possible to express $g$ as a product of fundamental reflections $\rgen_\alpha$, $\alpha\in \rootsys$ (equivalently, the length of the reduced Coxeter word representing $g$).  
Let $\inner{\cdot}{\cdot}$ be an inner product   of the Euclidean space in which the root system is embedded.
If $\lattice$ is a lattice of points invariant under a Weyl group $\group$ then the Weyl chamber (or fundamental region) of the group is the set of points 
$$
\chamber(\group)=\{x\in \lattice \,:\, \text{$\forall \alpha\in \Delta$,  $\inner{x}{\alpha}>0$} \}\,.
$$
The walls of the chamber $\chamber(\group)$ is the set of points 
$$\cwall(\group)=\{x\in \lattice \,:\, \text{$\forall \alpha\in \Delta$,  $\inner{x}{\alpha}=0$ and $\forall \beta\in \Delta\setminus\alpha$,  $\inner{x}{\beta}\ge 0$} \}.$$
In the affine case the chamber becomes an alcove.

The Gessel-Zeilberger theorem expresses the number of lattice paths (using steps from some step set) that are confined to the Weyl chamber can be expressed in terms of lattice paths  using steps from the same step set  that are \emph{not} confined to the chamber. 

Thus, let  $\cpath^{\group}_m(\ba\to \bb)$ be  the number of  length $m$ paths on $\lattice$ with step set $\stepset$ which start at $\ba$, end at $\bb$  and all   points of  are points the Weyl chamber $\chamber(\group)$. 
Let $\fpath_m(\ba\to \bb)$ be  the number of  length $m$ paths on $\lattice$ with step set $\stepset$ which start at $\ba$, end at $\bb$ (with no chamber constraint).

%
\begin{theorem}[Gessel-Zeilberger \cite{gessel:1992kx}]\label{thm_GZ}
Let $\lattice$ be a lattice of points, $\stepset$ the path step set and $\group$ a Weyl group satisfying
\begin{enumerate}
    \item[i)] $\lattice$  is invariant under $\group$, and
    \item[ii)] $\stepset$ is invariant under $\group$, and
    \item[iii)] For all $\alpha\in \Delta$ and for all $s\in\stepset$, $\inner{\alpha}{s}=\pm k(\alpha)$ where $k(\alpha)$ is a fixed number that only depends on $\alpha$.
\end{enumerate}
 is given by
\begin{equation}\label{eq_orbit}
   \cpath^{\group}_m(\ba\to \bb)
   =\sum_{g\in\group}(-1)^{\len(g)}\,\fpath_m(g(\ba)\to \bb)
\end{equation}
where $g(a)$ is the action $g$ on the lattice point $a$.
\end{theorem}
Note, the  theorem \emph{assumes} the Weyl chamber associated with a lattice path problem  is known. This paper gives a method for determining the Weyl chamber and hence whether or not the theorem can be used.
Since the step set is a subset of $\lattice$ if $\group$ leaves $\stepset$ invariant is also leaves $\lattice$ invariant.
Condition iii) enforces the constraint the every path contributing to $\fpath_m(g(a)\to b)$ that cross a wall of the chamber can only do so by stepping onto the wall first ie.\ it cannot ``step over'' a chamber wall.

Thus Theorem \ref{thm_GZ} can only give a solution to a short step quarter plane path problem if the quarter plane boundary corresponds to the walls of the Weyl chamber of the Weyl group of the step set (and lattice).

We will determine the Weyl group of the above four short step problems. 
Only the \rkre\ paths   have a Weyl chamber corresponding to the quarter plane. For the \kre\  paths the quarter plane us the union of two chambers (and the dividing wall). 
For the \gess\ paths the quarter plane is the union of   three adjacent chambers (and the  chamber walls).
For the \mr\ paths the group is affine and hence the quarter plane is the union of an infinite number of chambers.




\section{The step set and Weyl groups}

\newcommand{\action}{ \circ }
\newcommand{\zb}{\mathbf{0}}

The following   follows from elementary combinatorics.
  The free  \rkre, \kre\  \gess\ and \mr\ $t$-path   fixed length generating functions are
\begin{align}
G_t^{\rkm}(\ba;x,y)&=x^{a_1}y^{a_2}(x+y+\xb\yb)^t\\
      G_t^\text{K}(\ba;x,y)&=x^{a_1}y^{a_2}(\xb+\yb+xy)^t\\
        G_t^\text{G}(\ba;x,y)&=x^{a_1}y^{a_2}(x+\xb+\xb\yb+xy)^t\\
   G_t^\text{MR}(\ba;x,y)&=x^{a_1}y^{a_2}(\xb y + x \yb+xy)^t
\end{align}
where $\ba=(a_1,a_2)$.
 
\begin{figure}
    \centering
   \includegraphics[width=6cm]{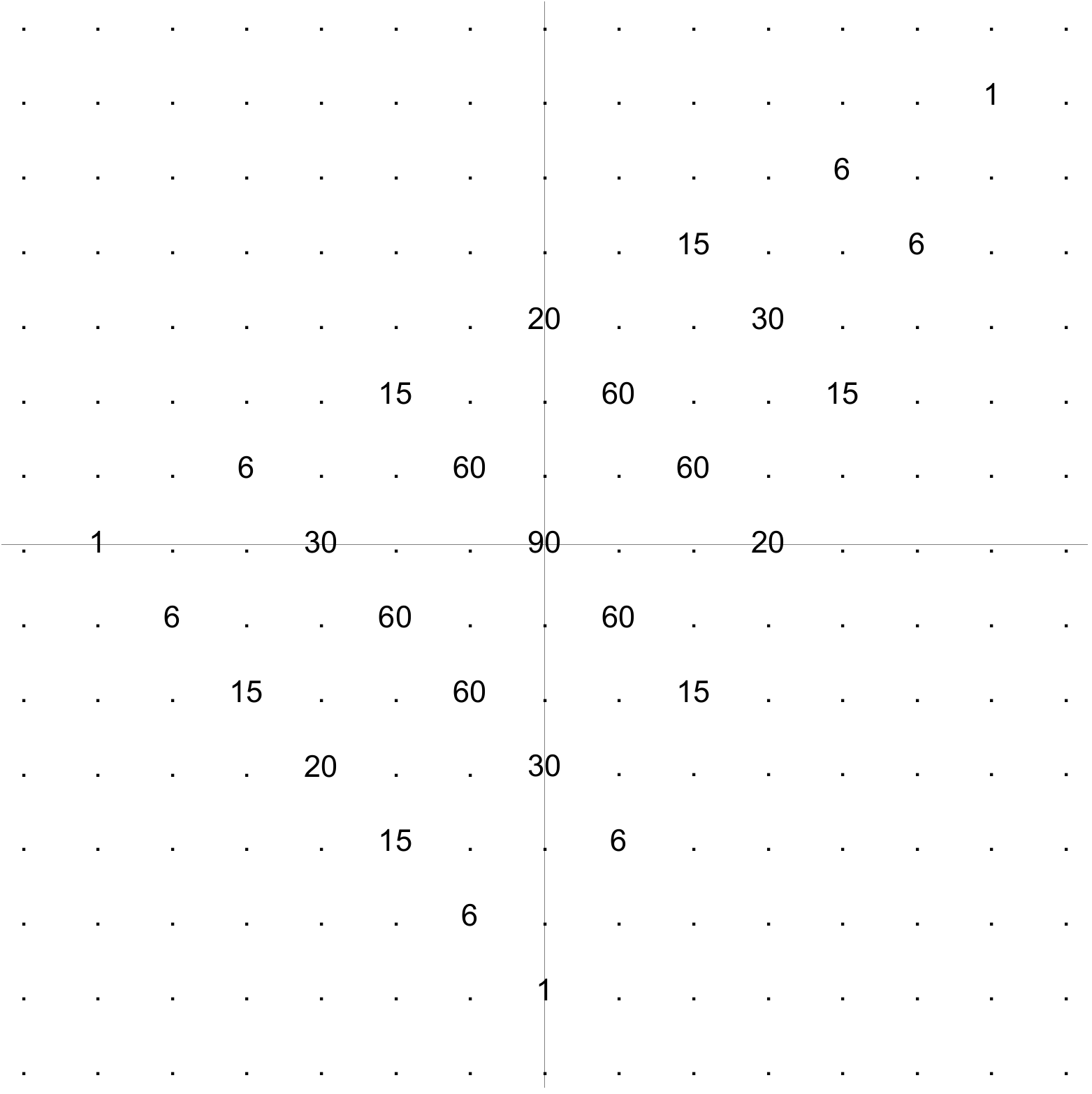}
      \includegraphics[width=6cm]{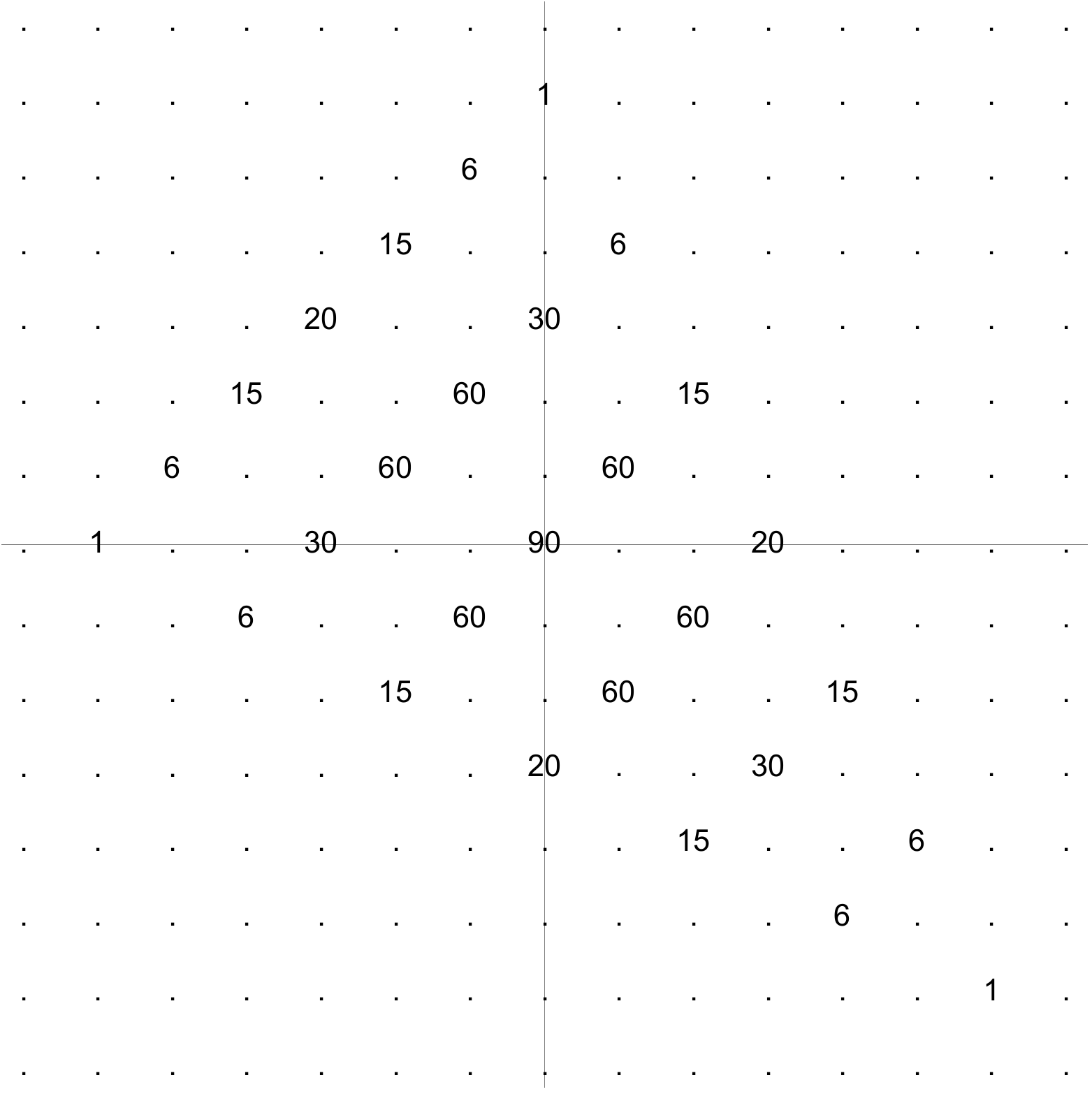}
    \caption{ Pascal slices for $t=6$ for (left) free \kre\ and (right) free \rkre\ paths. }
    \label{fig:free}
\end{figure}

The primary objects of interest are  $\gfreecount_t(\zb\to \bb)$, where $\zb=(0,0)$, which we will consider as defining integer valued counting functions on $\lattice$, that is
\begin{equation}\label{eq:freemapfunction}
      F_t^\text{P}\st \lattice\to\integer\,;\, (b_1,b_2)\mapsto   [x^{b_1} y^{b_2}]G_t^\text{P}(\zb;x,y) 
\end{equation}
The coefficients  of the free generating functions are trinomial coefficients giving the following  result.
\begin{prop}\label{prop_trin} The free  \rkre,  \kre\, \gess\ and \mr\ $t$-path  counting functions are:
\begin{subequations}\label{eq:trinom}
\begin{align}
    F_t^{\rkm} (a,b)=&
    \binom{t}{\frac{1}{3}(t-2a-b),\frac{1}{3}(t+a+2b),\frac{1}{3}(t+a-b)}
    \label{eq:trinom1} \\
    F_t^\text{K}(a,b)=&
    \binom{t}{\frac{1}{3}(t-2a+b),\frac{1}{3}(t+a-2b),\frac{1}{3}(t+a+b)}
    \label{trinom2}\\
    F_t^\text{G}(a,b)=&\sum_{k\ge 0}
    \binom{t}{k,k-a+b,-k+\frac{1}{2}(a-2b+t),-k+\frac{1}{2}(a+t)  }
    \label{eq:trinom4}\\
    F_t^\text{MR}(a,b)=&
    \binom{t}{\frac{1}{2}(t-a),\frac{1}{2}(t-b),\frac{1}{2}(a+b), }
    \label{eq:trinom3}
\end{align}
\end{subequations}
where the multinomials are assumed to vanish if  the value is not a positive integer.
\end{prop}
Examples of theses  functions are shown in figure \ref{fig:free}. These plots are the analogues of rows in Pascal triangle (but now correspond  to slices of a ``Pascal'' cone). We will thus call them \textbf{Pascal  slices}.

\begin{definition}
Let  $F_t^\text{P}$ be the map \eqref{eq:freemapfunction}  and $g\colon \lattice\to \lattice$ a distance preserving bijection that leaves $ F_t^\text{P}$ (or the Pascal plot) invariant. 
The set of such maps will be called the \textbf{reflection  group},  or isometry group, of the step set $\stepset_P$.   
\end{definition}

For the \kre, \rkre\  and \gess\ paths we shall see that the same finite group leaves $ F_t^P$ invariant \emph{independent} of $t$, whilst for \mr\ paths the group is different for every path length $t$ (but they are all subgroups of a single infinite group - an affine dihedral group).

If $G$ is a group of isometries of $\lattice$ then the action of $g\in G$   on $ F_t^\text{P}$, denoted  $g\act  F_t^\text{P}$ is defined in the usual way as
\begin{equation}
    (g\act  F_t^\text{P})(\bb) =  F_t^\text{P}(g^{-1}\act  \bb).
\end{equation}

Since $F_t^\text{P}(\bb)$ is given by trinomial coefficients the isometry group is related to the symmetry of the trinomial coefficients. Clearly the trinomial
\begin{equation}\label{eq:trinomdef}
    \binom{t}{n_1,n_2,n_3}=\frac{t!}{n_1!n_2!n_3!}
\end{equation}
is invariant under any permutation of $n_1$, $n_2$ and $n_3$. Thus to find the isometry groups of each of the   types of free path we look for isometries of $\lattice$ that permute the lower three parameters of the trinomials \eqref{eq:trinomdef}.  
For example, for the free \kre\ paths from \eqref{trinom2}:
\begin{subequations}
\begin{align}
    n_1(a,b)=&\frac{1}{3}(t-2a+b),\\
    n_2(a,b)=&\frac{1}{3}(t+a-2b),\\
    n_3(a,b)=&\frac{1}{3}(t+a+b)
\end{align}
\end{subequations}
and using the reflection
\[
g\colon (a,b)\mapsto (a,b)\begin{pmatrix} -1 &-1\\ 0 & 1
\end{pmatrix}
\]
give immediately the permutation
\[
n_1\mapsto n_1,\qquad n_2\mapsto n_3,\qquad n_3\mapsto n_2.
\]
and hence $g$ leaves \eqref{trinom2} unchanged.

In the case of the Gessel paths   it is only the sum of multinomials, \eqref{eq:trinom4}, that must be invariant. Thus a permutation of the trinomial indexes or a permutation of the trinomials in the sum would suffice. We only consider the former.

\subsection{Weyl groups for \rkre\ and \kre\ paths}

\paragraph{\rkre.}
The six possible permutations of the trinomial parameters give the following theorem.

\begin{theorem} \label{thm:krev}
The  isometry group of the free \rkre\ paths is an order six  dihedral group $D_3$ and have  the representation
\begin{align}
    g_0=&\begin{pmatrix} 1 &0\\ 0 & 1 \end{pmatrix},\quad  &g_1=\begin{pmatrix} -1 & 1\\ 0 & 1 \end{pmatrix} \qquad  g_2=&\begin{pmatrix}  1 &0\\ 1 & -1
\end{pmatrix}
\\
 g_3=&\begin{pmatrix} -1 & 1\\ -1 & 0 \end{pmatrix},\quad 
&g_4=\begin{pmatrix} 0 &-1\\ 1 & -1 \end{pmatrix} \qquad 
g_5=&\begin{pmatrix} 0 &-1\\ -1 & 0 \end{pmatrix} 
\end{align}
\end{theorem}
The elements $g_1$ and $g_2$  represent reflections with reflection lines having the Cartesian equations $x=0$ and $y=0$ respectively. As a Coxeter system $D_3$ is generated by the two reflections $g_1$ and $g_2$ and has presentation
\begin{equation}
     g_1^2=g_0\qquad g_1^2=g_0\qquad (g_1g_2)^3=g_0
\end{equation}
Thus $g_4=g_1g_2$ and $g_5=g_2g_1$ are order three rotations.

The two generating reflections are \emph{not} perpendicular to their reflection lines but reflect  at an oblique angle. 
The directions of the reflections are given by the left eigenvectors  corresponding to  eigenvalue $-1$ which are  $(-2,1)$ and $(-1,2)$ for $g_1$ and $g_2$  respectively. 

\newcommand{\bilin}[2]{\langle#1,#2\rangle}
\newcommand{\balpha}{{\boldsymbol{\alpha}}}

The oblique reflections can be seen as non-diagonal bilinear functional as follows. 
Conventionally reflections $s_\balpha$  across a hyperplane (a line in this case) defined by  a vector $\balpha$  corresponds to the map
\[
s_{\balpha}\colon \ba\mapsto \ba - \frac{2}{\bilin{\balpha}{\balpha} } \bilin{\ba}{\balpha} \balpha
\]
where the bilinear functional $\bilin{\cdot}{\cdot}$ is symmetric and positive definite and represented by the matrix $\mathbf{B}$  
\[
\bilin{\ba}{\bb}= \ba\cdot  \mathbf{B}\cdot  \bb^{T}
\]
where $ \bb^{T}$ is the transpose.
For the \rkre\ the matrix is
\[
\mathbf{B}^{\rkm} =\begin{pmatrix} 1 & -\frac{1}{2} \\
-\frac{1}{2} &  1
\end{pmatrix}\,.
\]
A basis for $\lattice$ can be chosen that would diagonalise $\mathbf{B}$ and hence change the reflections to perpendicular but this would then mean the quarter plane boundaries would not be perpendicular.

\begin{figure}
    \centering
 \includegraphics[width=12cm]{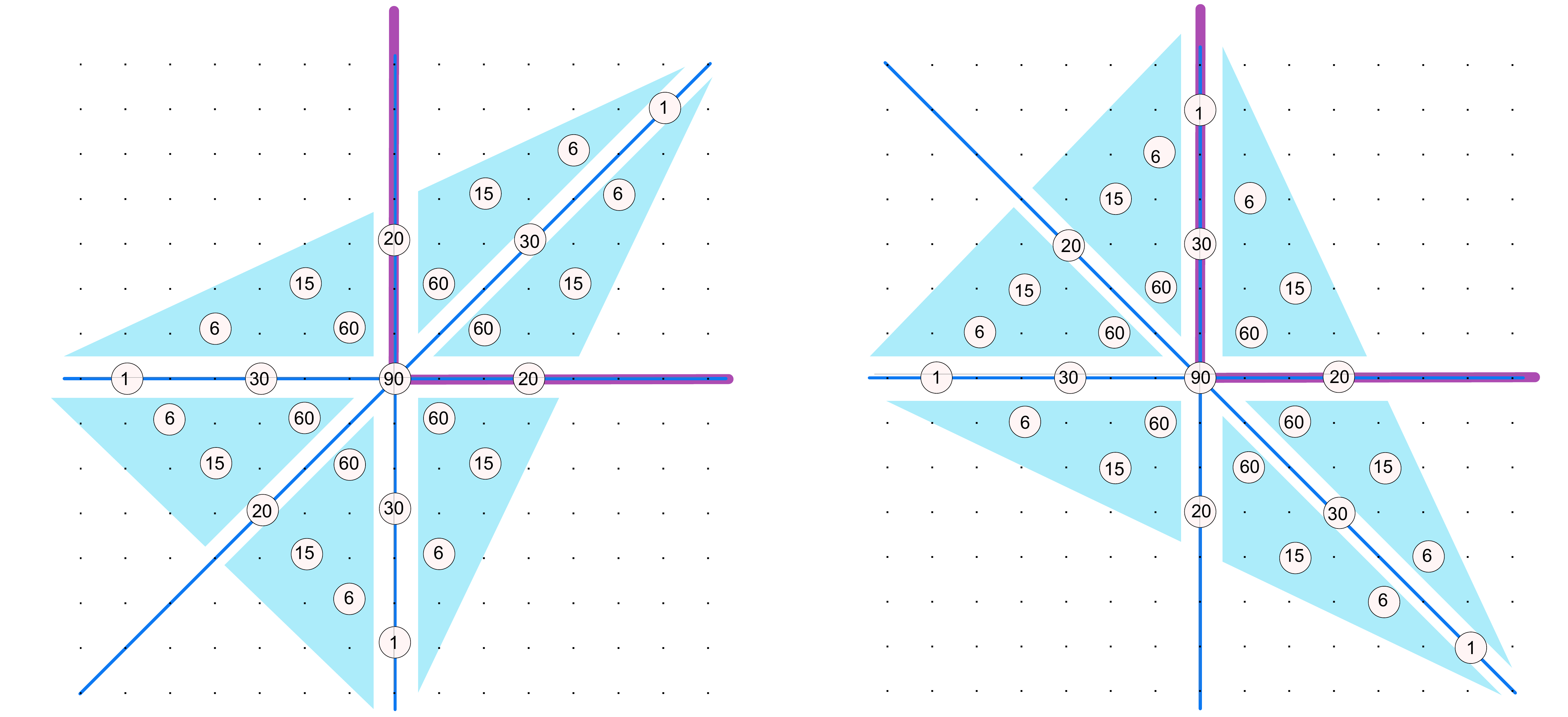}
    \caption{Weyl Chambers: (left) \kre\ and (right) \rkre}
    \label{fig:freeRKlines}
\end{figure}

The dihedral symmetry is readily seen in the Pascal slices - figure \ref{fig:freeRKlines} (right) shows the case for $t=6$ and shows the three reflection lines.
 
\paragraph{\kre.}
For   free \kre\   paths the result is similar, differing by a rotation of $\pi$ from the \kre\ counting function. The trinomial parameters are:
\begin{subequations}\label{eq:kretripara}
\begin{align}
    n_1(a,b)=&\frac{1}{3}(t-2a+b),\\
    n_2(a,b)=&\frac{1}{3}(t-2a+b),\\
    n_3(a,b)=&\frac{1}{3}(t+a+b)
\end{align}
\end{subequations}
Isometries of $\lattice$ which permute $n_1$, $n_2$ and $n_3$ are given by the following theorem.
\begin{theorem}\label{thmkrevvv}
The elements of the isometry group of the free \kre\ is an order six  dihedral group $D_3$ and have have representation
\begin{align}
g_0=&\begin{pmatrix} 1 &0\\ 0 & 1
\end{pmatrix},\qquad 
g_1=&\begin{pmatrix} -1 &-1\\ 0 & 1
\end{pmatrix} \qquad 
g_2=&\begin{pmatrix}  1 &0\\ -1 & -1
\end{pmatrix}\\
 g_3=&\begin{pmatrix} -1 &-1\\ 1 & 0
\end{pmatrix},\qquad 
g_4=&\begin{pmatrix} 0 & 1\\ -1 & -1
\end{pmatrix} \qquad 
g_5=&\begin{pmatrix} 0 & 1\\ 1 & 0
\end{pmatrix}\\
\end{align}
\end{theorem}
The proof of the theorem is a simple verification of the  permutations of \eqref{eq:kretripara}  using the matrices.

 The elements $g_1$ and $g_2$  represent  oblique  reflections with reflection lines having the Cartesian equations $x=0$ and $y=0$ respectively. The reflections are \emph{not} perpendicular to the lines but are at an oblique angle. The direction of the reflection is given by the eigenvector  corresponding to eigenvalue $-1$ which are  $(2,1)$ and $(1,2)$ for $g_1$ and $g_2$  respectively.
The dihedral symmetry is readily seen from a Pascal slice -   \figref{fig:freeRKlines} (left) shows the case for $t=6$ and shows the three reflection lines.

\paragraph{Weyl Chambers}

Rather than giving a formal definition of Weyl chambers (which can be found in any textbook on reflection/Coxeter groups) we give a brief informal description. For $\lattice$ the orbits of points in $\lattice$ under the action of  a dihedral group partitions $\lattice$ is subsets defined by the orbits. 
For $\lattice$ the orbits of size less than $n$   correspond to  points on any reflection line. The  set  of points in orbits of size $n$ partition into disjoint ``cones'' which are separated by the reflection lines. The action of the group on any cone will generate all the points of the other cones. 
Conventionally, a  cone whose points have positive coordinates is chosen - this is the Weyl chamber. 
For \rkre\ and \kre\ the cones and Weyl chamber can be seen in   \figref{fig:freeRKlines} and correspond to the set of points between the lines $y=x$ and $y=0$ and $x=0$ and $y=0$ respectively.

The essential difference between \rkre\ and \kre\  is apparent from figure \figref{fig:freeRKlines}. 
For the \rkre\ paths the quarter plane corresponds exactly to the Weyl chamber for the free paths, whilst for the \kre\ paths the quarter plane is bisected by a reflection line. 
This means a theorem such as that due to Gessel and Zielberger \cite{gessel:1992kx} cannot be used for the latter path problem (which requires the chamber walls correspond to the boundary of the domain).

It also means is that the standard Bethe Ansatz orbit Ansatz,
\begin{equation}\label{eq:bethe}
     G_t^{\rkm}(\ba\to \bb) =F_t^{\rkm}(\bb)   \sum_{g\in D_3} (-1)^{|g|} x^{g^{-1}(a_1)} y^{g^{-1}(a_2)} 
\end{equation}
where $\ba=(a_1,a_2)$, $D_3$ is given by the representation in \thmref{thm:krev} and $|g|$ is the number of reflections $g$, solves the \rkre\ enumeration problem (ie.\ gives the correct boundary conditions). 

In the case of \kre\ 
\begin{equation}
     G_t^{\rkm}(\ba\to \bb)\ne F_t^{\rkm}(\bb)   \sum_{g\in D_3} (-1)^{|g|} x^{g^{-1}(a_1)} y^{g^{-1}(a_2)} 
\end{equation}
where $D_3$ is given by the representation in \thmref{thmkrevvv},
the right hand side enforces a line of zeros along \emph{all} reflection lines. 
Since the \rkre\ has no reflection line through the quarter plane \eqref{eq:bethe} is the solution, however for the \kre\ there is reflection line bisecting the quarter plane and thus the orbit Ansatz does not solve \kre.  

This gives some insight as to why \kre\ is so much more difficult to  solve. The method used by \bm\ uses a ``half-orbit'' sum which almost solves the problem except for a ``line of defects''. 
This line can be removed by using Gessel factorisation \cite{bousquet-melou:2005uq}. 


\subsection{Weyl group for \gess\ paths}

A Pascal slice for $n=8$ is shown below. From this example it is clear there are four lines of reflection with oblique angles of reflection across all reflection line. 
In this form not all the angles between the lines of reflection are rational multiples of $\pi$, however the rotations are rational multiples of $\pi$ (see \figref{fig-gesselGrpRots}) as required of a finite dihedral group.

%
\begin{center} 
    \includegraphics[width=10cm]{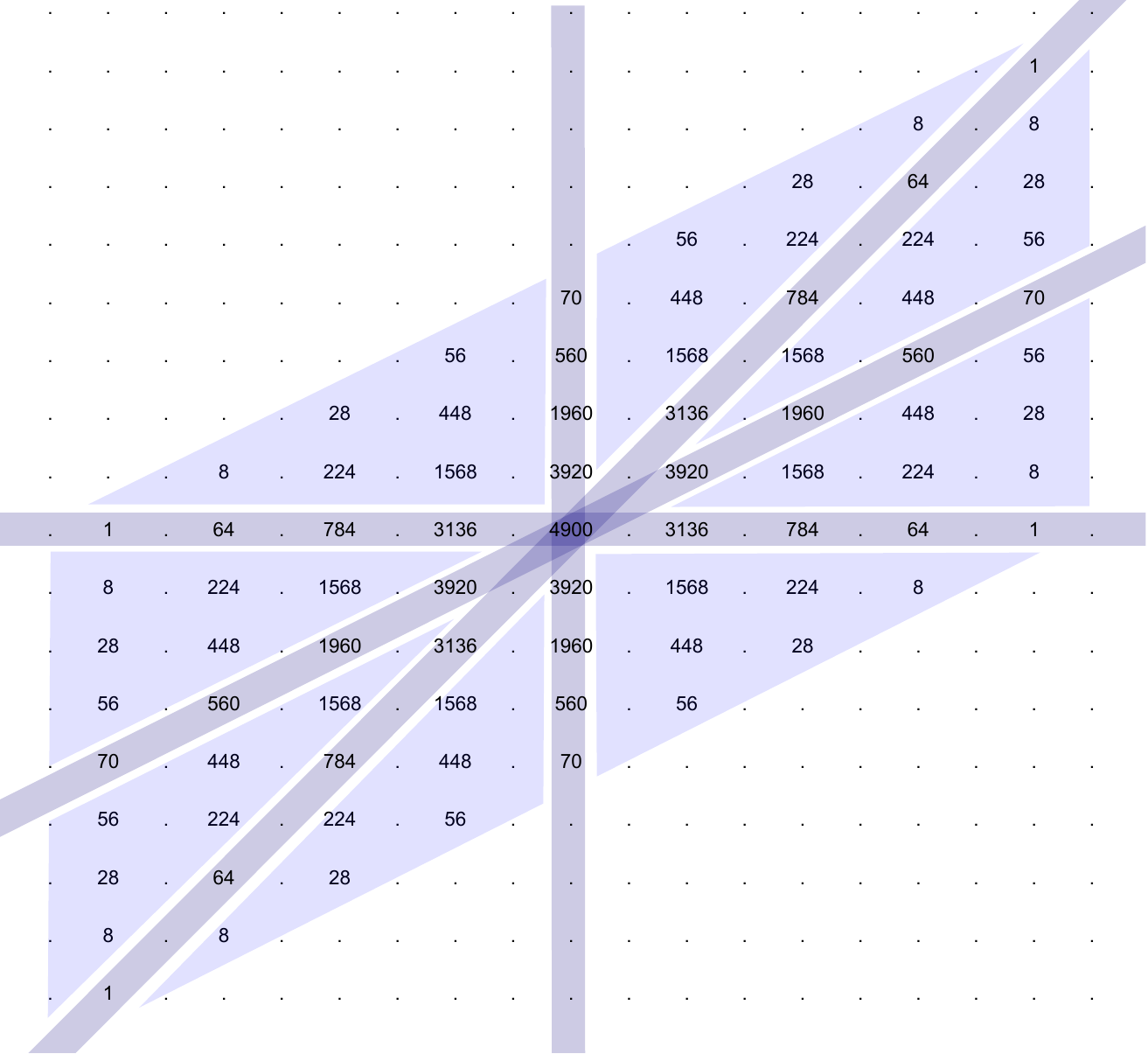}
\end{center}
The  isometry group of the free \gess\ paths is an order eight  dihedral group $D_4$ which have matrix representations  given in \figref{fig-gesselGrpRots} and \figref{fig:gesselGrpRefl}.
\begin{figure}[h]
    \centering
\begin{tabular}{|c|c|c|}
\hline 
$g_5=g_1g_4$ & $g_6=g_1g_3$ & $g_7=g_1g_2$\\
\hline  
& & \\
$  \begin{pmatrix} -1 &0 \\ 0 & -1
\end{pmatrix}$ 
& $  \begin{pmatrix}  1 & 1\\ -2 & -1
\end{pmatrix}$   & 
$  \begin{pmatrix} -1 &-1\\ 2 & 1
\end{pmatrix}$ \\
& & \\
$\displaystyle \frac{ \pi}{2}$  & 
$\displaystyle \frac{ \pi}{4} $ & 
$\displaystyle  \frac{ 5\pi}{4} $ \\
& & \\
\hline
\end{tabular}
    \caption{The three rotation elements of the  order eight  dihedral group $D_4$  of the free \gess\ paths  showing the representations  and the    angle of  rotation.}
    \label{fig-gesselGrpRots}
\end{figure}
\begin{figure}[th]
    \centering
\begin{tabular}{|c|c|c|c|}
\hline
$g_1$& $g_2$ & $g_3$ & $g_4$  \\
\hline
& & &  \\ 
$ \begin{pmatrix} 1 & 0\\ -2 & -1
\end{pmatrix}$ &  
$  \begin{pmatrix}  -1 &-1\\ 0 & 1
\end{pmatrix}$ & 
$ \begin{pmatrix} 1 &1\\ 0 & -1
\end{pmatrix}$   & 
$  \begin{pmatrix} -1 & 0\\ 2 & 1
\end{pmatrix} $  \\
& & & \\ 
$\displaystyle (1,1)\,, (1,0)$  & 
$\displaystyle (2,1)\,, (0,1)$ & 
$\displaystyle (0,1)\,, (2,1)$ & 
$\displaystyle (1,0)\,, (1,1)$ \\ 
& & & \\ 
$\displaystyle -1\qquad 1$  &
$\displaystyle -1\qquad 1$  &
$\displaystyle -1\qquad 1$  &
$\displaystyle -1\qquad 1$   \\
& & & \\ 
$y=0$& 
$x=0$ &
$y=x/2$  &
$y=x$ 
\\
& & & \\ 
$\displaystyle \frac{  \pi}{4}$ &
$\displaystyle \frac{\pi}{2}-\gamma$ &
$\displaystyle \gamma$ &
$\displaystyle  \frac{ 3 \pi}{4}$ \\
& & & \\
\hline
\end{tabular}
    \caption{The four reflection elements of the  order eight  dihedral group $D_4$  of the free \gess\ paths  showing the representations,  the  left  eigenvectors, the eigenvalues, the line of reflection  and the oblique  angle of  reflection where $\gamma=\arctan (2)= 1.10715\dots$}
    \label{fig:gesselGrpRefl}
\end{figure}

The Gessel path Weyl chamber is the set of points between the line $y=0$ and $y=x/2$.
Similar to the \kre\ paths,   the quarter plane does not coincide with the Weyl chamber so neither the Gessel-Zeilberger theorem nor the Bethe Ansatz will solve the problem.

\newpage

\subsection{Affine dihedral group for \mr\ paths}

The group structure of \mr\ paths is more complex. From \eqref{prop_trin}
\[
 F_t^\text{MR}(a,b)= 
    \binom{t}{\frac{1}{2}(t-a),\frac{1}{2}(t-b),\frac{1}{2}(a+b), } 
\]
The step set barycenter is at $(\frac{1}{3},\frac{1}{3})$ giving a path drift of $(\frac{t}{3},\frac{t}{3})$. 

From the trinomial,  $F_t^\text{MR}(a,b)$, is invariant if the isometry permutes the three functions
\begin{subequations}\label{eq:mrtripara}
\begin{align}
    n_1(a,b)=& t- a ,\\
    n_2(a,b)=& t - b ,\\
    n_3(a,b)=&  a+b \,.
\end{align}
\end{subequations}
\begin{figure}
    \centering
    \includegraphics[width=9cm]{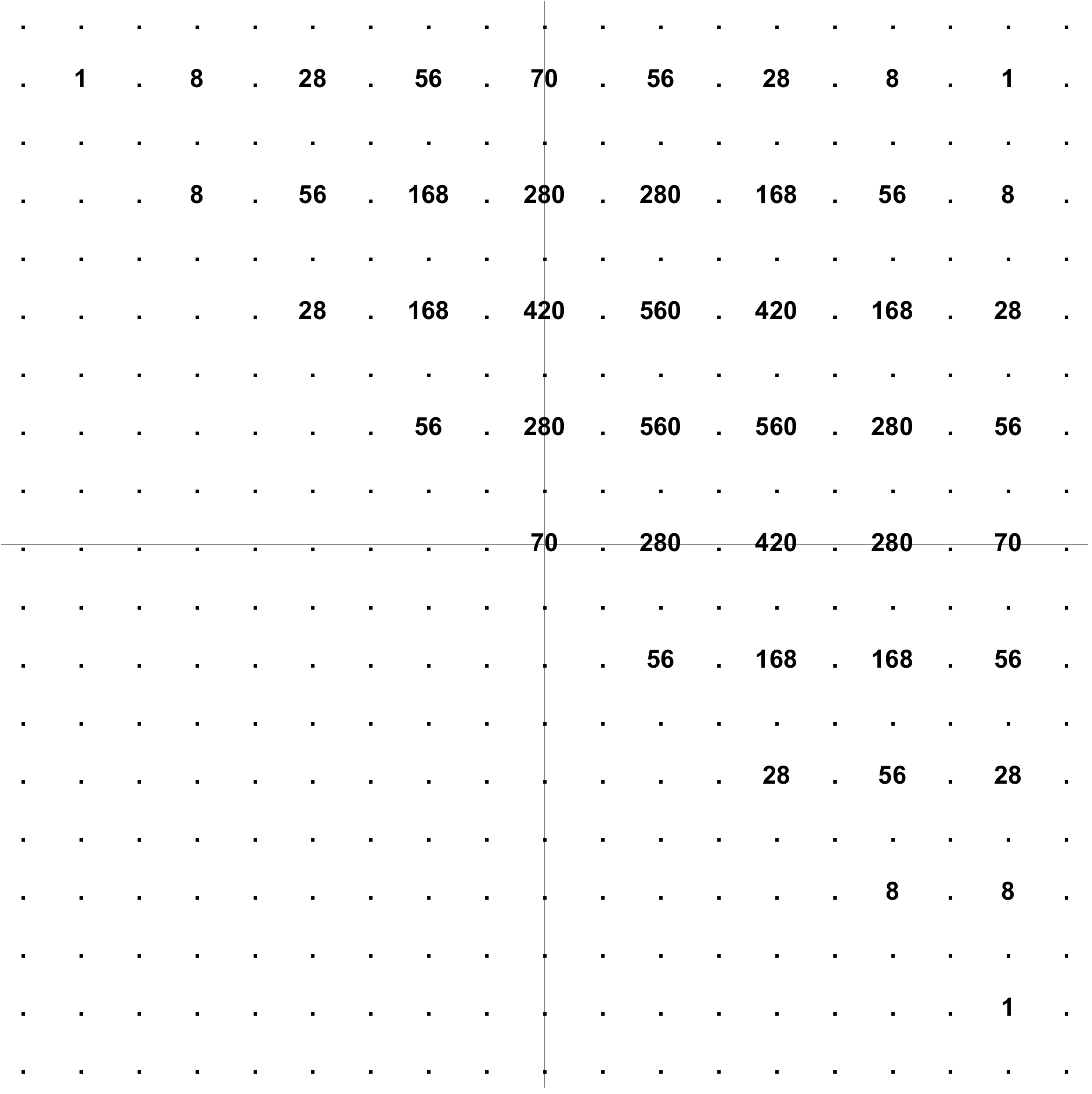}
    \caption{A Pascal slice for \mr\ paths corresponding to $t=8$.}
    \label{fig:mrpascalff}
\end{figure}
A Pascal slice is shown in \figref{fig:mrpascalff}.
These three equations differ fundamentally from the previous path problems in that one equation, $n_3$, does not depend on $t$. 
The absence of $t$ from $n_3$ forces us to extend the isometries from reflections and rotations to include translations, that is, to use affine maps. For example, if we try
\[
g_{321}\colon (a,b)\mapsto (a,b)\begin{pmatrix} 1 &-1\\ 0 & -1
\end{pmatrix}+  t(\tau_1,\tau_2)
\]
where $(\tau_1,\tau_2)$ is some row vector. Substituting into \eqref{eq:mrtripara} we see that the three equations are permuted 
\[
 n_1\mapsto n_3\,,\quad   
 n_2\mapsto n_2\quad  \text{and}\quad  
 n_3\mapsto n_1
\]
only  if $\tau_1=0$ and $\tau_2=1$. Using affine isometries gives the following theorem.
\begin{theorem}
The elements of the isometry group of the free \mr\ paths is an affine  dihedral group $\check{D}_3$ where every element of the group is of the form
\begin{equation}
     T_\sigma(k_1,k_2)\colon (a,b)\mapsto (a,b)g_\sigma+(k_1,k_2)
\end{equation}
where the linear maps $g_\sigma$ have representations:
\begin{align}
g_{123}=&\begin{pmatrix}  1 & 0\\ 0 &  1
\end{pmatrix},\qquad 
g_{132}=&\begin{pmatrix}  1 & -1\\ 0 & -1
\end{pmatrix} \qquad 
g_{321}=&\begin{pmatrix} -1 & 0\\ -1 &  1
\end{pmatrix}\\
 g_{312}=&\begin{pmatrix} -1 &  1\\ -1 & 0
\end{pmatrix},\qquad 
g_{231}=&\begin{pmatrix}  0 & -1\\  1 & -1
\end{pmatrix} \qquad 
g_{213}=&\begin{pmatrix} 0 &1\\ 1 & 0
\end{pmatrix}\,.
\end{align}
The \mr\ trinomial parameters \eqref{eq:mrtripara} are permuted by the subgroup of elements
\begin{align}
     T_{123}(0,0)& &
     T_{213}(0,0)& &
     T_{132}(0,t)\notag\\
     T_{321}(t,0)& &
     T_{312}(t,0)&  &
     T_{231}(0,t)
\end{align}
\end{theorem}
Again, the proof of the theorem is a simple substitution of the affine maps into \eqref{eq:mrtripara} and showing each permutes $n_1$, $n_2$ and $n_3$. 
Since the group is affine we now have Weyl alcoves rather than chambers. Since the alcoves are a finite set of elements they cannot correspond to the quarter plane.
\begin{figure}
    \centering
    \includegraphics[width=9cm]{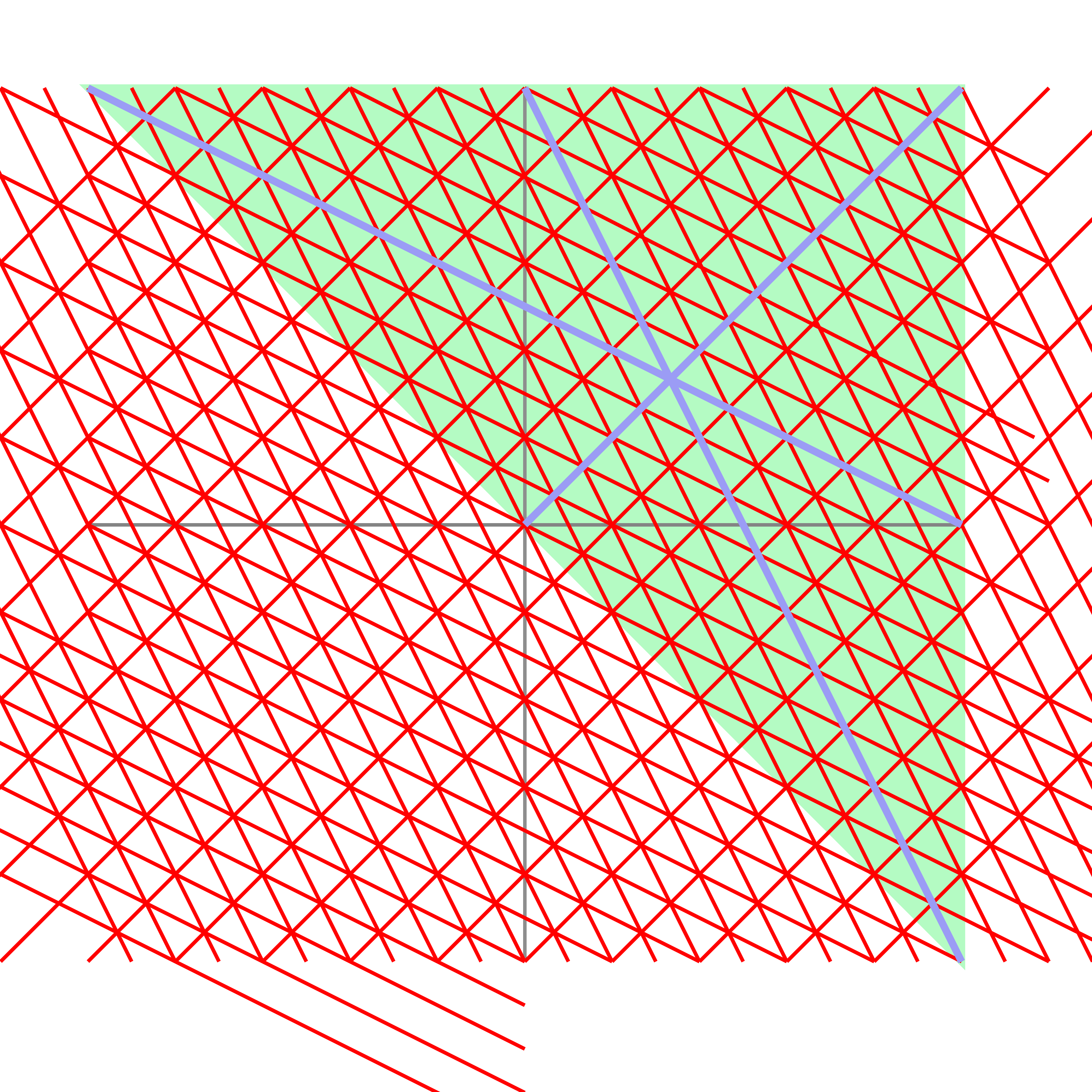}
    \caption{The alcoves of the affine $\check{D}_3$ group of the \mr\ paths. The blue lines are the reflection planes which leave the $t=10$ (hence they intersect at $(\frac{10}{3},\frac{10}{3})$)  Pascal slice invariant.}
    \label{fig:mrpascal}
\end{figure}

It is not clear how this group structure can be used in conjunction with some sort of orbit Ansatz to solve the enumeration problem.

%
\section{Acknowledgement} We would like to   thank the Australian Research Council (ARC) and the Centre of Excellence for Mathematics and Statistics of Complex Systems (MASCOS) for financial support. 

\newpage
\bibliographystyle{unsrt}
\bibliography{ThreePaths}


\end{document}